\documentclass[sigconf]{acmart}

\usepackage{enumerate}
\usepackage{color}
\usepackage{url}
\usepackage{subfigure}
\usepackage{graphicx,epstopdf}
\usepackage{float}
\usepackage{amsmath,  amssymb, algorithmic, algorithm, url, color, bbm, enumitem}
\usepackage[maxfloats=30]{morefloats}

\newcommand{\hide}[1]{}
\newcommand{\diffd}{\delta}

\usepackage{ifthen}
\newboolean{showcomments}
\setboolean{showcomments}{true}
\newcommand{\yue}[1]{\ifthenelse{\boolean{showcomments}}
{ \textcolor{red}{(Yue says:  #1)}}{}}
\newcommand{\josh}[1]{\ifthenelse{\boolean{showcomments}}
{ \textcolor{red}{(Josh says:  #1)}}{}}
\newcommand{\zhenhua}[1]{\ifthenelse{\boolean{showcomments}}
{ \textcolor{red}{(Zhenhua says:  #1)}}{}}
\newcommand{\addcites}[0]{\ifthenelse{\boolean{showcomments}}
{ \textcolor{green}{(add citation(s))}}{}}
\newcommand{\addcite}[0]{\ifthenelse{\boolean{showcomments}}
{ \textcolor{green}{(add citation(s))}}{}}
\newcommand{\addref}[0]{\ifthenelse{\boolean{showcomments}}
{ \textcolor{green}{(add ref)}}{}}
\newcommand{\todo}[1]{\ifthenelse{\boolean{showcomments}}
{ \textcolor{red}{(To do:  #1)}}{}}
\newcommand{\fixes}[1]{\ifthenelse{\boolean{showcomments}}
{\textcolor{black}{#1}}{}}

  \makeatletter

\graphicspath{{images/}}

\usepackage{booktabs} % For formal tables

\begin{document}
	
\copyrightyear{2017}
\acmYear{2017}
\setcopyright{rightsretained}
\acmConference{SIGMETRICS '17}{}{June 5--9, 2017, Urbana-Champaign, IL, USA}
\acmPrice{}
\acmDOI{http://dx.doi.org/10.1145/3078505.3078546}
\acmISBN{ACM ISBN 978-1-4503-5032-7/17/06}

%Conference
%\copyrightyear{2017}
%\acmYear{2017}
%\setcopyright{rightsretained}
%\acmConference{SIGMETRICS '17}{June 05-09, 2017}{Urbana-Champaign, IL, USA}\acmDOI{http://dx.doi.org/10.1145/3078505.3078546}
%\acmISBN{978-1-4503-5032-7/17/06}
%\acmPrice{}

\title{Incentivizing Reliable Demand Response with Customers' Uncertainties and Capacity Planning}
%\titlenote{Produces the permission block, and copyright information}
%\subtitle{Submission \#7}
%\subtitlenote{The full version of the author's guide is available as \texttt{acmart.pdf} document}

%\author{Submission \#7}
%\affiliation{\institution{Affiliation}}
%\email{email}

\author{Joshua Comden}
\affiliation{\institution{Stony Brook University}}
\email{joshua.comden@stonybrook.edu}

\author{Zhenhua Liu}
\affiliation{\institution{Stony Brook University}}
\email{zhenhua.liu@stonybrook.edu}

\author{Yue Zhao}
\affiliation{\institution{Stony Brook University}}
\email{yue.zhao.2@stonybrook.edu}

\maketitle

\vspace{-0.1in}
\section{Introduction}

%%%%%%%%%%%%%%%%%%%%%%%%%%%%%
% Problems with uncertianty %
%%%%%%%%%%%%%%%%%%%%%%%%%%%%%

One of the major issues with the integration of renewable energy sources into the power grid is the increased uncertainty and variability that they bring.
%~\cite{sgdoe}. 
%The limited capability to accurately predict this variability makes it challenging for the load serving entities (LSEs) to respond to it.
%~\cite{nrelvariability}. 
If this uncertainty
%variability
is not sufficiently addressed, it will limit the further penetration of renewables into the grid and even result in blackouts.
Compared to energy storage, Demand Response (DR) has advantages to provide reserves to the load serving entities (LSEs) in a cost-effective and environmentally friendly way.
DR programs work by changing customers' loads when the power grid experiences a contingency such as a mismatch between supply and demand.
%The decision that must be made by DR programs is how much load demand should each customer change.
Uncertainties from both the customer-side and LSE-side make designing algorithms for DR a major challenge.
%LSEs make predictions about the net load demand and purchase capacity to dispatch controllable supply accordingly, but they do not know the true mismatch between supply and demand until the time arrives.
%On the other hand, customers who are accustomed to having electricity supply on demand are not be able to accurately estimate how much disutility a future change in load demand would bring them.
%These uncertainties mean that the LSE does not know how exactly how much aggregate DR and reserve capacity it will need, and the customers do not know how much change they will be willing to provide.
%For this reason, we propose two simple DR contracts between the customers and the LSE.

This paper makes the following main contributions:
%This paper first models the problem of joint capacity planning and demand response program design with a stochastic optimization problem.
%which incorporates the uncertainties from renewable energy generation, customer power demands, as well as the customers' costs in providing DR.
(i) We propose
%online
DR control policies based on the optimal structures of the offline solution.
(ii) A distributed algorithm is developed for implementing the control policies without efficiency loss.
(iii) We further offer an enhanced policy design by allowing flexibilities into the commitment level. 
(iv) We perform real world trace based numerical simulations
%Results
which demonstrate that the proposed algorithms can achieve near optimal social cost.
Details can be found in our extended version \cite{comden2017incentivizingarXiv}.
\vspace{-0.05in}
%\section{Model}
\section{Optimization Problem}
\label{sec:model}

The goal is %therefore
to \emph{simultaneously} decide the capacity planning $\kappa$ and a practical DR policy $\mathbf{x}(D,\boldsymbol{\delta})$ to minimize the \emph{expected} social cost caused by a random aggregate supply-demand mismatch $D$ (which captures mismatches from both the generation side and the load side). 
%\eqref{eq:expected-cost}.
%\begin{subequations}\label{opt:general}
%	\begin{align}
%		\min_{\kappa,\mathbf{x}(\delta_i,\delta_r)} & C_\text{cap}(\kappa)\nonumber \\
%		&+\mathbb{E}_{\delta_i,\delta_r,C_i(\cdot)}\left[\sum_{i}C_i(x_i(\delta_i,\delta_r))+C_\text{g}\left(D-\sum_{i}x_i(\delta_i,\delta_r)\right)\right] \nonumber \\
%		\text{s.t. } & \max_{\delta_i,\delta_r}\left\{D-\sum_{i}x_i(\delta_i,\delta_r)\right\}\leq \kappa \label{const:CAPpolicy1} \\
%		&\min_{\delta_i,\delta_r}\left\{D-\sum_{i}x_i(\delta_i,\delta_r)\right\}\geq -\kappa. \label{const:CAPpolicy2}
%	\end{align}
%\end{subequations}
\vspace{-0.05in}
\begin{subequations}\label{opt:general}
	\begin{align}
		\min_{\kappa,\mathbf{x}(D,\boldsymbol{\delta})} & C_\text{cap}(\kappa)\nonumber \\
		&+\mathbb{E}_{D,\boldsymbol{\delta},\mathbf{C}(\cdot)}\left[\sum_{i}C_i(x_i(D,\delta_i))+C_\text{g}\left(D-\sum_{i}x_i(D,\delta_i)\right)\right] \nonumber \\
		\text{s.t. } & \max_{D,\boldsymbol{\delta}}\left\{D-\sum_{i}x_i(D,\delta_i)\right\}\leq \kappa \label{const:CAPpolicy1} \\
		&\min_{D,\boldsymbol{\delta}}\left\{D-\sum_{i}x_i(D,\delta_i)\right\}\geq -\kappa. \label{const:CAPpolicy2}
	\end{align}
\end{subequations}
where $\delta_i$ and $C_i(\cdot)$ are respectively for customer $i$ the individual random demand mismatch and random cost function (e.g. $a_ix_i^2$ with $a_i$ as a random coefficient) for performing DR, $C_\text{cap}(\cdot)$  %\yue{is this correct?}
and $C_\text{g}(\cdot)$ are respectively the LSE's cost for purchasing capacity and for managing the remaining mismatch,
We note that \eqref{const:CAPpolicy1} and \eqref{const:CAPpolicy2} are worst-case constraints so that the remaining mismatch does not go beyond the purchased capacity.
%To optimize the policy is known to be challenging, so we first provide the upper and lower bound to this optimization.
The two main challenges of this problem are (i) deciding the optimal capacity $\kappa$ \emph{before} implementing the DR policy, and (ii) optimizing an online DR policy. %instead of a fixed amount.
The cost functions are assumed to be convex. %We make the following mild assumption:
\vspace{-0.05in}
\subsection*{Optimal Real-time Solution} %\yue{this is not a policy, but a solution}}

\begin{figure}
	\subfigure[LIN]{{\includegraphics[width=0.40\columnwidth]{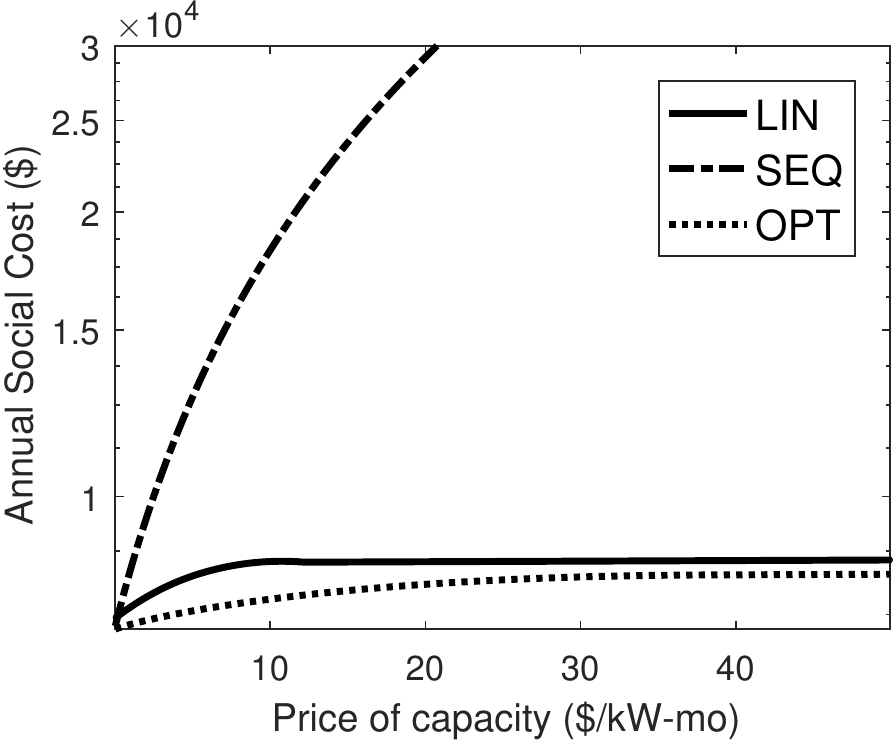}}
		\label{FIG:COMPSCH_SOC}}
	\subfigure[LIN$^+(\rho)$]{{\includegraphics[width=0.40\columnwidth]{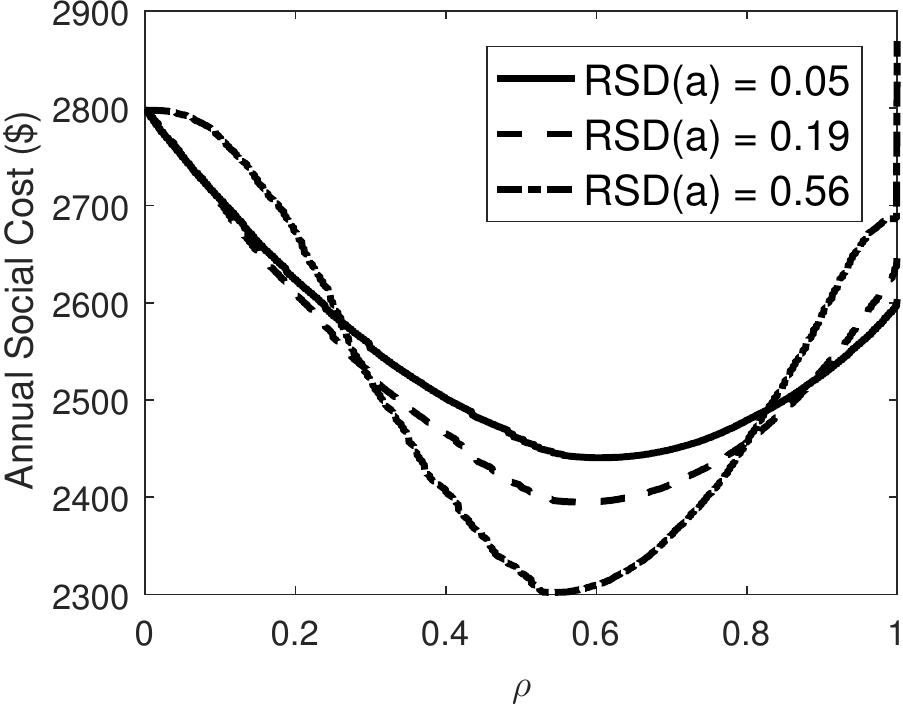}}
		\label{FIG:FDR_SOC}}
	\vspace{-0.2in}
	\caption{Annual Social Cost vs. (a) price of capacity in LIN compared to different baselines, (b) level of commitment in LIN$^+(\rho)$ for different amounts of Relative Standard Deviations (RSD) on the customer cost parameter $a$.}
	\vspace{-0.15in}
\end{figure}

\vspace{-0.05in}
We provide the characterization of the optimal real-time solution to reveal special structures that we take advantage of in our policy design (Section \ref{sec:alg}).
%The first key result regarding the problem is the convexity, as stated formally in Theorem~\ref{th:offline-convexity}.
%The convexity is crucial for our proposed algorithm in Section~\ref{sec:alg}.
%\begin{theorem}
%\label{th:offline-convexity}
%	\eqref{opt:OFFLINE1} is a convex optimization problem over $\kappa$. %\yue{over $\kappa$ to be exact?}
%\end{theorem}
%\begin{proof}
%	The result follows from Assumption \ref{as:kappaconvex} and Lemma \ref{lm:convexKAPPA} given below.
%\end{proof}
%The proof requires the following lemmas and we restate the real-time decision (i.e. inside the expectation) of Problem \eqref{opt:OFFLINE1} as:
%\zhenhua{Josh, please update the following proof for the convexity of \eqref{opt:OFFLINE1}.}
The real-time DR decision problem for a given capacity $\kappa$ at a time $t$ is:
\vspace{-0.15in}
\begin{subequations}\label{opt:OFFLINE1_REALTIME}
	\begin{align}
		R(\kappa;t):=\min_{\mathbf{x}(t)} \quad & \sum_{i}C_i(x_i(t);t)+C_\text{g}\left(D(t)-\sum_{i}x_i(t)\right) \\
		\text{s.t.} \quad & -\kappa\leq D(t)-\sum_{i}x_i(t)\leq\kappa. \label{const:Capacity_RT}
	\end{align}
\end{subequations}
\vspace{-0.2in}
\begin{lemma}\label{thm:convexity2}
	Problem \eqref{opt:OFFLINE1_REALTIME} is a convex optimization problem.
	%, \eqref{opt:GENall2D}, \eqref{opt:LSEalg2}, and \eqref{opt:CUSTalg2} are convex optimization problems.
\end{lemma}
\vspace{-0.05in}

The Karush-Kuhn-Tucker (KKT) optimality conditions of this real-time problem show that
when the capacity constraint %\eqref{eq:KKTgen2prim}
on $\kappa$ is non-binding, i.e., $-\kappa<D(t)-\sum_{j\in\mathcal{V}}x_j^*(t)<\kappa$,
it implies that $C'_i(x_i^*(t))=C'_\text{g}(D(t)-\sum_{j\in\mathcal{V}}x_j^*(t))$.
This means that the marginal cost for each customer to provide demand response is the same, all of which is equal to the LSE's marginal cost to tolerate the mismatch.
Furthermore, we get the following lemma which helps determine the optimal capacity in the next subsection:
\vspace{-0.05in}
\begin{lemma}\label{lm:convexKAPPA}
	$R(\kappa;t)$ as defined by \eqref{opt:OFFLINE1_REALTIME} is a convex function of $\kappa$. Additionally the negative of the sum of dual variables $\underline{\theta}+\overline{\theta}$ from constraint \eqref{const:Capacity_RT} %that results from $(\kappa;t)$ 
	is the subgradient of $R(\kappa;t)$ w.r.t. $\kappa$.
\end{lemma}

\vspace{-0.15in}
\subsection*{Optimal Capacity}
\vspace{-0.05in}
%A lower bound on the minimum social cost is given by the \emph{offline/a-posteriori} optimal solution. Specifically, the offline optimum is given by the following:
We can use the real-time decision problem \eqref{opt:OFFLINE1_REALTIME} to decide what the optimal capacity should be in the following capacity problem:
%\begin{subequations}\label{opt:OFFLINE1}
%	\begin{align}
%		\min_{\kappa}\text{ } & C_k(\kappa) + \mathbb{E}_{\delta_i,\delta_r}\min_{\mathbf{x}(t)}\left\{\sum_{i}C_i(x_i(t);t)+C_\text{g}\left(D(t)-\sum_{i}x_i(t)\right)\right\} \nonumber \\
%		\text{s.t. } & -\kappa\leq D(t)-\sum_{i}x_i(t)\leq\kappa,\quad \text{for each realization}~t. \label{const:Capacity}
%	\end{align}
%\end{subequations}
\begin{align}
	\min_{\kappa} \quad C_\text{cap}(\kappa) + \mathbb{E}_{D,\boldsymbol{\delta},\mathbf{C}(\cdot)}\left[R(\kappa;t)\right] \label{opt:OFFLINE1} %\\
	%\text{s.t. } & -\kappa\leq D(t)-\sum_{i}x_i(t)\leq\kappa,\quad \text{for each realization}~t. \label{const:Capacity}
\end{align}
\vspace{-0.05in}
%Note that, the minimization over $\mathbf{x}(t)$ is performed \emph{inside} the expectation, meaning that it is performed after observing the realizations of the random variables. This results in the fact that the offline optimum can never be beaten by any online policy $\mathbf{x}(\delta_i, \delta_r)$. 
\begin{theorem}
	\label{th:offline-convexity}
	\eqref{opt:OFFLINE1} is a convex optimization problem over $\kappa$. %\yue{over $\kappa$ to be exact?}
\end{theorem}
\vspace{-0.05in}
%\zhenhua{Joshua, please update the conditions for $\kappa^*$}

The KKT optimality conditions for the capacity problem and Lemma \ref{lm:convexKAPPA} give us the following result:
%\vspace{-0.05in}
\begin{equation}\label{eq:POSToptKappa}
	C'_\text{cap}(\kappa^*) = \mathbb{E}_{D,\boldsymbol{\delta},\mathbf{C}(\cdot)}\left[\theta(\kappa^*;t)\right]
\end{equation}
where we use the notation of $\theta(\kappa;t)$ as a function to represent the sum of the optimal dual variables for constraint \eqref{const:Capacity_RT}.
%$\underline{\theta}^*+\overline{\theta}^*$,
%from the KKT conditions \eqref{eq:KKTgen2} for given $(\kappa;t)$.
This means that for an optimal capacity, the marginal cost of capacity must equal the expected dual price for that capacity constraint.

\section{Policy Design}
\label{sec:alg}

\subsection*{Linear policy}
\label{sec:alg_lin}
\vspace{-0.05in}
Motivated by the desire to find a simple DR policy $\mathbf{x}(D,\delta_{i\in\mathcal{V}})$ that preserves convexity, we focus on a simple but powerful linear demand response policy that is a function of total and local net demands:
%\begin{align}
%	x_i({\diffd}_i)=\alpha_iD+\beta_i{\diffd}_i+\gamma_i. \label{pol:GENlinear}
%\end{align}
\vspace{-0.15in}
\begin{align}
	x_i(D,{\diffd}_i)=\alpha_iD+\beta_i{\diffd}_i+\gamma_i. \label{pol:GENlinear}
\end{align}
%which makes 
%\todo{need to change ${\diffd}_i$ to the error, not absolute}
Intuitively, there are three components: $\alpha_iD$ implies each customer shares some (predefined) fraction of the global mismatch $D$; $\beta_i{\diffd}_i$ means  customer $i$ may need to take additional responsibility for the mismatch due to his own demand fluctuation and estimation error; finally, $\gamma_i$, the constant part, can help when the random variables $\mathbb{E}[D]$ and/or $\mathbb{E}[\delta_i]$ is nonzero.
%Accordingly, constraint \eqref{const:CAPpolicy} becomes:
%
%\yue{The following can be made a lot simpler, also now that we have a lower bound in addition to upper bound, we need to update the following.}
%\begin{subequations}\label{const:KAPPAs}
%	\begin{equation}
%		\sum_{i\in\mathcal{V}}\kappa_i+\kappa_r\leq \kappa_\text{tot}
%	\end{equation}
%	\begin{equation}
%		\underline{\diffd}_i\left(1-\sum_{j\in\mathcal{V}}\alpha_j-\beta_i\right)-\gamma_i\leq \kappa_i,\quad \forall i\in\mathcal{V}
%	\end{equation}
%	\begin{equation}
%		\overline{\diffd}_i\left(1-\sum_{j\in\mathcal{V}}\alpha_j-\beta_i\right)-\gamma_i\leq \kappa_i,\quad \forall i\in\mathcal{V}
%	\end{equation}
%	\begin{equation}
%		\underline{r}\left(\sum_{j\in\mathcal{V}}\alpha_j-1\right)\leq \kappa_r
%	\end{equation}
%	\begin{equation}
%		\overline{r}\left(\sum_{j\in\mathcal{V}}\alpha_j-1\right)\leq \kappa_r
%	\end{equation}
%\end{subequations}
Then the LSE needs to solve \eqref{opt:general} with \eqref{pol:GENlinear} to obtain the optimal parameters for the linear contract, i.e., $\boldsymbol{\alpha}, \boldsymbol{\beta}, \boldsymbol{\gamma}$, as well as the optimal capacity $\kappa$.
%Then the LSE needs to solve the following optimization problem to obtain the optimal parameters for the linear contract, i.e., $\boldsymbol{\alpha}, \boldsymbol{\beta}, \boldsymbol{\gamma}$, as well as the optimal capacity $\kappa$:
%\begin{align}
%	\min_{\boldsymbol{\alpha},\boldsymbol{\beta},\boldsymbol{\gamma},\kappa}\text{ } &
%	C_\text{cap}\left(\kappa\right)+\sum_{i\in\mathcal{V}}\mathbb{E}_{\delta_i, \delta_r, C_i}\left[C_i(\alpha_iD+\beta_i{\diffd}_i+\gamma_i)\right]\nonumber \\
%	 & \quad +\mathbb{E}_{\delta_i, \delta_r}\left[C_\text{g}\left(D-\sum_{i\in\mathcal{V}}(\alpha_iD+\beta_i{\diffd}_i+\gamma_i)\right)\right] \label{opt:GENall2} \\
%	\text{s.t. } & \eqref{const:CAPpolicy1}, \eqref{const:CAPpolicy2}\nonumber
%\end{align}

%\todo{again change the above $d$ to error.}

%The expectation can be separated by customer and LSE because of the linearity of expectation property.

%\todo{Add the theorem to say \eqref{opt:GENall2} is convex and the proof. }
\vspace{-0.05in}
\begin{theorem}
	Problem \eqref{opt:general} with the linear policy \eqref{pol:GENlinear}
	%\eqref{opt:GENall2}
	is a convex optimization problem.
\end{theorem}
\vspace{-0.2in}
\subsection*{Distributed algorithm}
\vspace{-0.05in}

In most cases, the LSE's information on the customers' cost functions is much less accurate than the customer themselves'. This can also be due to privacy concerns. %and cannot solve \eqref{opt:general} with \eqref{pol:GENlinear}
%\eqref{opt:GENall2}
%for capacity planning and the parameters in the linear contract.
To handle this, we design a distributed algorithm
%to decompose the problem,
so that the LSE does not need the information of the customer cost functions, %each customer does not have to reveal her private cost function 
while still achieving
%is solving an optimization problem based on her own cost function, and the LSE is solving another optimization problem without knowing the customers' cost function. 
%The goal is to achieve
the optimal $(\kappa^{LIN}, \boldsymbol{\alpha}^*,\boldsymbol{\beta}^*,\boldsymbol{\gamma}^*)$ for Problem \eqref{opt:general} with the linear policy \eqref{pol:GENlinear}.
%\yue{LIN is not defined}.
%that can be used in LIN.
We introduce and substitute $(u_i,v_i,w_i)$ for $(\alpha_i,\beta_i,\gamma_i)$ in each customer's estimated cost function $\hat{C}_i(\cdot)$
and the LSE uses the corresponding price set $(\pi_i,\lambda_i,\mu_i)$ to incentivize each customer to change their parameters.

%\vspace{0.05in}
\noindent
\fbox{\begin{minipage}{26em}
		\textbf{Distributed Algorithm for LIN}:
		\begin{enumerate}[start=0]
			\item \textbf{Initialization:} $(\boldsymbol{\alpha},\boldsymbol{\beta},\boldsymbol{\gamma},\mathbf{u},\mathbf{v},\mathbf{w},\boldsymbol{\pi},\boldsymbol{\lambda},\boldsymbol{\mu}):=\mathbf{0}$.
			\item \textbf{LSE:} receives $(u_i,v_i,w_i)$ from each customer $i\in\mathcal{V}$.
			\begin{itemize}
				\item Solves Problem \eqref{opt:LSEalg2}	and updates $(\boldsymbol{\alpha},\boldsymbol{\beta},\boldsymbol{\gamma})$ with the optimal solution.
				\item Updates the stepsize:
				\vspace{-0.05in}
				\begin{equation}\label{step:STEPupdate}
				\eta=\frac{\zeta/k}{||(\boldsymbol{\alpha},\boldsymbol{\beta},\boldsymbol{\gamma})-(\mathbf{u},\mathbf{v},\mathbf{w})||_2}
				\end{equation}
				where $\zeta$ is a small constant and $k$ is the iteration number.
				\item Updates the dual prices, $\forall i\in\mathcal{V}$:
			\end{itemize}
			\begin{small}
				\begin{equation}\label{step:DUALupdate}
				(\pi_i,\lambda_i,\mu_i):=(\pi_i,\lambda_i,\mu_i)+\eta\left((\alpha_i,\beta_i,\gamma_i)-(u_i,v_i,w_i)\right)
				\end{equation}
			\end{small}
			\begin{itemize}
				\item Sends $(\pi_i,\lambda_i,\mu_i)$ to the each customer respectively.
			\end{itemize}
			\item \textbf{Customer $i\in\mathcal{V}$:} receives $(\pi_i,\lambda_i,\mu_i)$ from LSE.
			\begin{itemize}
				\item Solves Problem \eqref{opt:CUSTalg2} and updates $(u_i,v_i,w_i)$ with optimal solution.
				\item Sends $(u_i,v_i,w_i)$ to the LSE.
			\end{itemize}
			\item Repeat Steps 1-2 until $||(\boldsymbol{\alpha},\boldsymbol{\beta},\boldsymbol{\gamma})-(\mathbf{u},\mathbf{v},\mathbf{w})||_2\leq \epsilon$ where $\epsilon$ is the tolerance on magnitude of the subgradient.
		\end{enumerate}
	\end{minipage}}
Thus $\pi_iu_i+\lambda_iv_i+\mu_iw_i$ is the total payment to customer $i$ for the linear demand response policy. %a constant reduction and response to the overall state of the system. 
%Accordingly, \eqref{opt:GENall2} is decomposed as
The individual customer's %optimization
problem for a given set of prices is
\begin{align}
	\min_{u_i,v_i,w_i} \mathbb{E}_{D,\diffd_i}\left[\hat{C}_i(u_iD+v_i{\diffd}_i+w_i)\right]-\pi_iu_i-\lambda_iv_i-\mu_iw_i \label{opt:CUSTalg2}
\end{align}
\vspace{-0.05in}
while the LSE's optimization problem among all the customers is
\begin{align}
	\min_{\boldsymbol{\alpha},\boldsymbol{\beta},\boldsymbol{\gamma},\kappa}\text{ } & 	C_\text{cap}\left(\kappa\right)+\sum_{i\in\mathcal{V}}(\pi_i\alpha_i+\lambda_i\beta_i+\mu_i\gamma_i)\nonumber \\
	 & \quad  +\mathbb{E}_{D,\boldsymbol{\delta}}\left[C_\text{g}\left(\sum_{i\in\mathcal{V}}({\diffd}_i-\alpha_iD-\beta_i{\diffd}_i-\gamma_i)-r\right)\right] \label{opt:LSEalg2} \\
	 \text{s.t. } & \eqref{const:CAPpolicy1}, \eqref{const:CAPpolicy2}\nonumber
\end{align}
%\vspace{-0.05in}
%\begin{equation}
%\min_{\boldsymbol{\alpha},\boldsymbol{\beta},\boldsymbol{\gamma},\kappa} 	C_\text{cap}\left(\kappa\right)+\sum_{i\in\mathcal{V}}(\pi_i\alpha_i+\lambda_i\beta_i+\mu_i\gamma_i)\nonumber
%\end{equation}
%\begin{equation}\label{opt:LSEalg2}
%+\mathbb{E}_{\mathbf{\diffd},r}\left[C_\text{g}\left(\sum_{i\in\mathcal{V}}({\diffd}_i-\alpha_iD-\beta_i{\diffd}_i-\gamma_i)-r\right)\right]
%\end{equation}
%\begin{equation}
%\text{s.t. }\eqref{const:CAPpolicy1}, \eqref{const:CAPpolicy2}\nonumber
%\end{equation}
%Problems \eqref{opt:CUSTalg2} and \eqref{opt:LSEalg2} can be solved with standard stochastic optimization techniques such as the Stochastic Subgradient Method with Monte Carlo sampling \cite{boyd2014stochsubgradient}. With linear or quadratic cost functions, they in fact can be solved as a deterministic convex optimization problem whose parameters are determined by first and second order moments of $\mathbf{\diffd}$ and $\delta_r$.   
%To solve the decomposed problems, we must ensure the customers' and LSE's decisions satisfy \eqref{const:SEPARATE}.
In order for the customers and the LSE to negotiate and obtain the optimal prices we use
%We achieve this using
the Subgradient Method (see \cite{bertsekas1999nonlinear} Chapter 6).
%as follows:
%to obtain the optimal dual prices in the following
\vspace{-0.05in}
\begin{theorem}\label{th:convergence}
	%Given Assumption \ref{as:DUALbound},
	%and \ref{as:SUBGRADbound},
	The distributed algorithm's %best
	trajectory of dual prices
	%\yue{confusing, what does best dual prices mean?}
	converge to the optimal dual prices for Problem \eqref{opt:general} with \eqref{pol:GENlinear}.
\end{theorem}

\vspace{-0.15in}
\subsection*{Flexible Commitment Demand Response}
\label{sec:fcdr}
\vspace{-0.05in}

%\subsubsection{LIN with flexibility: LIN$^+(\rho)$}

One potential drawback of LIN is that customers are forced to follow the specified linear policy.
In some cases, customers may face a very high cost to follow the policy, e.g., when there are some critical jobs to be finished, represented by a larger $a_i(t)$.
%From the social perspective, it is not beneficial to force a customer under much higher than average costs to provide demand response, as this would force them to incur higher costs compared to the LSE tolerating more mismatch.
Motivated by this observation and some existing regulation service programs,
%described in Section~\ref{sec:bg},
we modify the LIN policy to add some flexibility limited by a single parameter $\rho$.
We call the new algorithm LIN$^+(\rho)$
%Under LIN$^+(\rho)$,
where each customer has up to $1-\rho$ (in percentage) of the time slots in which they do not need to follow the policy according to her realized $\alpha_i(t)$.
In other words, she may let 
% pick their original $d_i(t)$ with
$x_i(t)=0$ for such timeslots.
%Another caveat is that we only allow the customers to pick such timeslots according to their cost functions, e.g., through the parameters $\alpha_i(t)$, instead of the real-time $D(t)$.
%The reason is that from the social perspective, allowing the customers with higher $\alpha_i(t)$ to violate results in lower social costs, but a higher $D(t)$ is seen by all customers and this is when the LSE needs DR the most.
%Allowing some to violate means others need to take more, so there is a negative externality here.
%From classical economic theory, this may lead to significant efficiency loss, or formally, the price of anarchy.
%In reality, this can be enforced by calculating $\mathbb{E}_i[D|i\;\text{violates}]$ for the timeslots that a particular customer chooses to violate.
%If this value is significantly higher than $\mathbb{E}[D]$, then the customer faces some penalty.
Note that although we add the flexibility to LIN in this paper, the approach is in fact general and can be applied to a wide range of fully committed programs.
%as follows.
%Each customer is allowed to violate her commitment by up to $1-\rho$ (in percentage) of the timeslots.
%A customer can run a local optimization to decide during what timeslots to violate, while LSE can add some constraints to align the local optimization with the social optimization, such as those described above.

%\yue{Move the distributed algorithm above flexible DR}. 
%\vspace{-0.1in}
\section{Performance Evaluation}
\label{sec:linear-evaluation}
\vspace{-0.05in}

\subsubsection*{Experimental Setup}

%We aim to use realistic parameters in the experimental setup to evaluate the performance of our proposed algorithms, and to understand the properties of solutions to each algorithm.
We simulate an LSE supplying power to 300 customers.
%and a demand response timeslot that is five minutes long.
%The LSE must first purchase capacity for which we model the cost as a linear function $c\kappa$.
%The generation cost function for the LSE is modeled as a quadratic function. %$A\left(D-\sum_ix_i\right)^2$.
Each customer has a particular demand of load which we model by utilizing the traces obtained from the UMass Trace Repository% which give very granular load measurements from three homes
~\cite{barker2012smart}.

\vspace{-0.05in}
%\subsubsection*{How well do PRED and LIN perform?}
\subsubsection*{LIN is close to optimal} %\yue{remove PRED in the graph}}

%\todo{remove the original LIN and replace by LIN-TRN, we only need one line}

%Now we move to this key question of the evaluation. 
%To evaluate the benefit
%in terms of social cost savings for our algorithms, Figure~\ref{f.COMPSCH}
Figure \ref{FIG:COMPSCH_SOC} compares the social cost of LIN to baselines using the offline optimal OPT \eqref{opt:OFFLINE1} as a lower bound and sequential algorithm SEQ as an upper bound.
%The baseline OPT requires knowledge about each realization of the parameters, i.e., the exact renewable generation, customers' demand and cost functions at each timeslot, and therefore not practically feasible.
The baseline SEQ first makes a conservative capacity planning decision about $\kappa$, and then sets a price for DR to obtain a targeted amount of DR.
%Recall that
%We simply use this as a lower bound of the social cost.
%On the other hand, SEQ makes conservative capacity planning decision about $\kappa$ first, and then performs similar to PRED.
%By comparing to the cost of SEQ, the benefit of joint optimization of capacity planning on $\kappa$ and real-time demand response is highlighted.
The social cost of
%PRED and
LIN
%, shown in Figure~\ref{FIG:COMPSCH_SOC}, are
is no more than $10\%$ higher compared to the fundamental limit OPT
%\zhenhua{How much cost savings compared to SEQ?}
and is significantly less than SEQ.
The social cost of SEQ increases rapidly with increasing capacity prices because of the conservative 90kW capacity used by SEQ to protect the system from \emph{any} leftover mismatch.

\vspace{-0.05in}
%\subsubsection{Cost versus mismatch leftover tradeoff}
\subsubsection*{Additional cost savings brought by LIN$^+(\rho)$}
\label{sec:fcdr-evaluation}
%\todo{add FCDR to Figure~\ref{f.COMPSCH}}.

%Now we evaluate the additional cost savings brought by LIN$^+(\rho)$.
%Figure~\ref{f.FDR} highlights the tradeoff between cost and mismatch leftover that cannot be balanced due to capacity constraints.
%Recall that there is little mismatch leftover in LIN, or equivalently, LIN$^+$(1).

Depicted in Figure~\ref{FIG:FDR_SOC}, as $\rho$ decreases from 1, the social cost first decreases due to the fact that some customers with very high $a_i(t)$ are allowed to not provide demand response.
As $\rho$ continues to decrease, we have more customers not providing demand response and the cost actually goes up again.
This is because the LSE's penalty for the mismatch becomes larger than the costs of customers to provide demand response.
%\zhenhua{need to explain what the normalized numbers mean.}
%On the other hand, Figure~\ref{FIG:FDR_VF} highlights that the mismatch leftover increases with lower $\rho$ as more customers are allowed to not provide demand response.
%Importantly, the increase is very slow at the beginning when we decrease $\rho$ from $1$ because there is enough capacity to handle the small deficit in demand response when $1-\rho$ is small.
%However, the social cost decreases a lot.
%This actually highlights the great potential to decrease $\rho$ appropriately to achieve a lower social cost but still have little mismatch leftover.
%For instance, in our case study shown in Figure~\ref{f.FDR},
At $\rho=0.8$, it achieves a cost savings 7-8\%.
%with less than 1\% of mismatch leftover.
Recall that the gap between LIN and the offline optimal OPT is about 10\%. 
This means LIN$^+(\rho^*)$ achieves near optimal cost.

%Key messages:
%\begin{itemize}
%\item With lower $\rho$, social cost first decreases and then increases.
%\item Remaining mismatch increases with lower $\rho$, but very slowly at the beginning
%\item It is possible to pick some $\rho<1$ to achieve a good tradeoff.
%\end{itemize}

%\begin{figure}
%	\includegraphics[width=0.95\columnwidth]{FDR_soc_frq}
%	\caption{Normalized social cost and capacity violation frequency versus the level of commitment.}
%	\label{f.FDR}
%\end{figure}

%\input{S7-conclusion}
\vspace{-0.05in}
\section{Acknowledgments}
\vspace{-0.05in}
This research is supported by NSF grants CNS-1464388 and CNS-1617698.
\vspace{-0.05in}

\vspace{-0.05in}
\bibliographystyle{ACM-Reference-Format}
\bibliography{reference}


\begin{thebibliography}{00}

%%% the Web address.

\ifx \showCODEN    \undefined \def \showCODEN     #1{\unskip}     \fi
\ifx \showDOI      \undefined \def \showDOI       #1{{\tt DOI:}\penalty0{#1}\ }
  \fi
\ifx \showISBNx    \undefined \def \showISBNx     #1{\unskip}     \fi
\ifx \showISBNxiii \undefined \def \showISBNxiii  #1{\unskip}     \fi
\ifx \showISSN     \undefined \def \showISSN      #1{\unskip}     \fi
\ifx \showLCCN     \undefined \def \showLCCN      #1{\unskip}     \fi
\ifx \shownote     \undefined \def \shownote      #1{#1}          \fi
\ifx \showarticletitle \undefined \def \showarticletitle #1{#1}   \fi
\ifx \showURL      \undefined \def \showURL       #1{#1}          \fi
% invisible to TeX
\providecommand\bibfield[2]{#2}
\providecommand\bibinfo[2]{#2}
\providecommand\natexlab[1]{#1}
\providecommand\showeprint[2][]{arXiv:#2}

\bibitem[\protect\citeauthoryear{Barker, Mishra, Irwin, Cecchet, Shenoy, and
  Albrecht}{Barker et~al\mbox{.}}{2012}]%
        {barker2012smart}
\bibfield{author}{\bibinfo{person}{Sean Barker}, \bibinfo{person}{Aditya
  Mishra}, \bibinfo{person}{David Irwin}, \bibinfo{person}{Emmanuel Cecchet},
  \bibinfo{person}{Prashant Shenoy}, {and} \bibinfo{person}{Jeannie Albrecht}.}
  \bibinfo{year}{2012}\natexlab{}.
\newblock \showarticletitle{Smart*: An open data set and tools for enabling
  research in sustainable homes}.
\newblock \bibinfo{journal}{{\em SustKDD, August\/}}  \bibinfo{volume}{111}
  (\bibinfo{year}{2012}), \bibinfo{pages}{112}.
\newblock


\bibitem[\protect\citeauthoryear{Bertsekas}{Bertsekas}{1999}]%
        {bertsekas1999nonlinear}
\bibfield{author}{\bibinfo{person}{Dimitri~P Bertsekas}.}
  \bibinfo{year}{1999}\natexlab{}.
\newblock \bibinfo{booktitle}{{\em Nonlinear programming}}.
\newblock \bibinfo{publisher}{Athena scientific Belmont}.
\newblock


\bibitem[\protect\citeauthoryear{Comden, Liu, and Zhao}{Comden
  et~al\mbox{.}}{2017}]%
        {comden2017incentivizingarXiv}
\bibfield{author}{\bibinfo{person}{Joshua Comden}, \bibinfo{person}{Zhenhua
  Liu}, {and} \bibinfo{person}{Yue Zhao}.} \bibinfo{year}{2017}\natexlab{}.
\newblock \showarticletitle{Harnessing Flexible and Reliable Demand Response
  Under Customer Uncertainties}.
\newblock \bibinfo{journal}{{\em arXiv preprint arXiv:1704.04537\/}}
  (\bibinfo{year}{2017}).
\newblock


\end{thebibliography}

%%%%There is slack to expand before the "Distributed Algorithm for LIN" block.
%%%%There is slack to expand Section 6.3 "Flexibe Commitment Demand Response".

%\appendix
%\vspace{-0.1in}
%\section*{Appendix}
%\input{appendix}

\end{document}